\documentclass[reqno, 10pt]{amsart}

\usepackage{article-preamble}
\usepackage{amsmath,amssymb}
\hypersetup{colorlinks=true,
    linkcolor=blue,
    citecolor=blue,
    pdftitle={A decomposition theorem for balanced measures},
    pdfauthor={Gregory Baimetov, Ryan Bushling, Ansel Goh, Raymond Guo, Owen Jacobs, and Sean Lee}}
\begin{document}

\title{A decomposition theorem for balanced measures}
\author[Baimetov, Bushling, Goh, Guo, Jacobs, \and Lee]{Gregory Baimetov, Ryan Bushling, Ansel Goh, \\ Raymond Guo, Owen Jacobs, \and Sean Lee}

\address{Department of Mathematics \\ University of Washington, Box 354350 \\ Seattle, WA 98195-4350}
\email{baimetov@uw.edu}
\email{reb28@uw.edu}
\email{anselgoh@uw.edu}
\email{rpg360@uw.edu}
\email{obj3@uw.edu}
\email{seanl6@uw.edu}

\begin{abstract}
    Let $G = (V,E)$ be a finite, connected, simple graph. A probability measure $\mu$ on $V$ is called \textit{balanced} if it has the following property: if $T_\mu(v)$ denotes the ``earth mover's" cost of transporting all the mass of $\mu$ from all over the graph to the vertex $v$, then $T_\mu$ attains its global maximum at each point in the support of $\mu$. We prove a decomposition result that characterizes balanced measures as convex combinations of suitable ``extremal" balanced measures that we call \textit{basic}. An upper bound on the number of basic balanced measures on $G$ follows, and an example shows that this estimate is essentially sharp.
\end{abstract}

\maketitle

\section{Introduction} \label{s:intro}

\subsection{Balanced Measures} \label{ss:balanced}

Let $G = (V,E)$ be a finite, connected, simple graph and $\mu$ a probability measure on the vertices $V$, which we identify with a non-negative function $\mu \!: V \to \R_{\geq 0}$ satisfying $\sum_{v \in V} \mu(v) = 1$. If $d(u,v)$ denotes the distance between vertices $u$ and $v$, then it is natural to define the cost of transporting the $\mu$-mass at $u$ over to $v$ by $d(u,v) \mu(u)$. The \define{transport cost function} of $\mu$ is defined as the cost of relocating the \textit{entire} mass of $\mu$ to a single vertex
\begin{equation*}
    T_\mu(v) := \sum_{u \in V} d(u,v) \mu(u).
\end{equation*}
Equivalently, this is the Wasserstein-$1$ distance (or ``earth mover's distance") of $\mu$ from the unit point mass at $v$. The behavior of $T_\mu$ reflects both the distribution of $\mu$ and the global geometry of $G$, and studying those probability measures on $G$ with ``exceptional" transport cost functions sheds light on that geometric structure. To that end, $\mu$ is said to be a \define{balanced measure} on $G$ if the maximal transport cost is attained at every point in the support of $\mu$, i.e., if
\begin{equation} \label{eq:balanced-defn}
    \mu(v) > 0 \quad \text{implies} \quad T_\mu(v) = \max_{u \in V} T_\mu(u).
\end{equation}

\begin{center}
\begin{figure}[h]
    \begin{tikzpicture}
        \node at (0,0) {\includegraphics[width=0.45\textwidth]{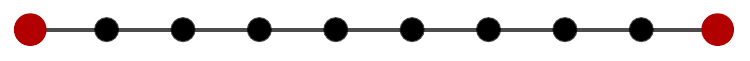}};
        \node at (-3,0.5) {$1/2$};
        \node at (3,0.5) {$1/2$};
    \end{tikzpicture}
    \vspace{-5pt}
    \caption{A balanced measure on the $n$-vertex path graph $P_n$. While a graph may admit many balanced measures, the equidistributed probability measure supported on the endpoints is in fact the unique balanced measure on the path graph.}
    \label{fig:path-graph}
\end{figure}
\end{center}

Letting $\spt \mu := \{ v \in V \!: \mu(v) > 0 \}$ be the support of $\mu$ and
\begin{equation*}
    M_\mu := \left\{ v \in V \!: T_\mu(v) = \max_{u \in V} T_\mu(u) \right\}
\end{equation*}
the \define{cost-maximizing set} of $\mu$, we equivalently require that $\spt \mu \subseteq M_\mu$. The reverse inclusion need not hold, and the relationship between the supports and cost-maximizing sets of balanced measures will be of major consequence in \S \ref{s:basic-measures}.

To illustrate the idea of a balanced measure, first consider the path graph $P_n$ in Figure \ref{fig:path-graph} and the measure $\mu$ that assigns equal mass to the endpoints. The cost of moving all of $\mu$ to a single vertex is independent of the destination, since $k/2 + (n-k)/2 = n/2$. Hence, the global transport cost is maximized at every point of $G$ and, in particular, at the endpoints of $P_n$ where $\mu$ is supported. As \eqref{eq:balanced-defn} holds trivially, $\mu$ is balanced. Steinerberger introduced balanced measures in \cite{steinerberger2023sums} and proved the following existence result.

\begin{thm}[Steinerberger \cite{steinerberger2023sums}] \label{thm:existence}
    Every connected, finite simple graph admits a balanced measure.
\end{thm}

\begin{center}
\begin{figure}[h!]
    \begin{tikzpicture}
        \node at (0,0) {\includegraphics[width=0.3\textwidth]{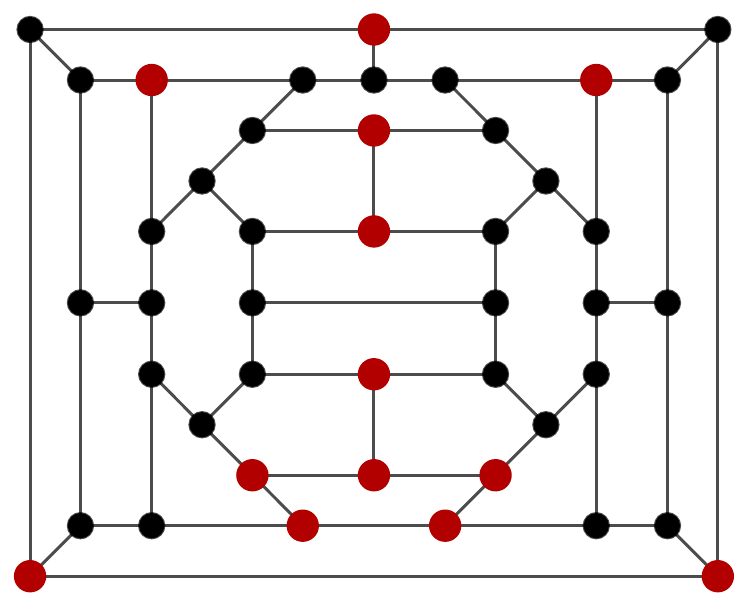}};
        \node at (6,0) {\includegraphics[width=0.3\textwidth]{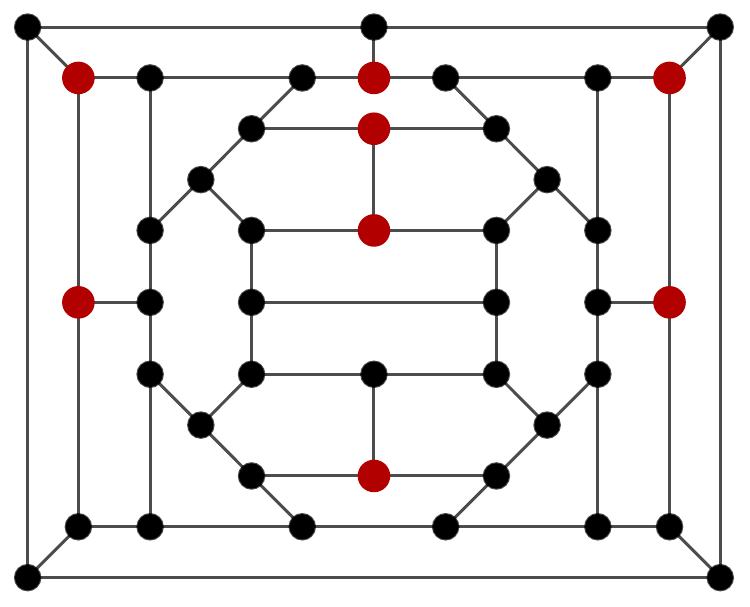}};
    \end{tikzpicture}
    \vspace{-5pt}
    \caption{The supports of two balanced measures on the $44$-vertex Grinberg graph.}
    \label{fig:grin}
\end{figure}
\end{center}

In fact, the result says a little more: if we start with an arbitrary finite list of vertices and then perform a greedy construction\textemdash adding the vertex that maximizes the sum of distances to the existing vertices\textemdash then the resulting sequence of empirical probability measures contains a subsequence converging to a balanced measure (cf.~\S \ref{ss:motivation}). Figure \ref{fig:grin} displays the supports of two balanced measures on the Grinberg graph obtained via this procedure. Among the natural questions that arise in \cite{steinerberger2023sums} are the size and structure of the set of all balanced measures on a graph, and the ease with which these measures can be generated algorithmically.

\vs{-0.75}

\subsection{A Structure Theorem}

While \cite{steinerberger2023sums} guarantees existence of balanced measures, it says little about their structure. The main motivation behind our paper is to shed light on this structure by relating the different balanced measures on a graph to each other. We begin with some terminology. Recall that $\spt \mu$ is the support of a measure and $M_\mu$ is its cost-maximizing set.
\begin{defn}
    We say that a probability measure $\mu$ on the vertices of a graph $G$ is \define{basic} if every other probability measure $\nu \neq \mu$ has different support or different maximizing set, i.e., if at least one of $\spt \mu \neq \spt \nu$ or $M_\mu \neq M_\nu$ holds.
\end{defn}
As an example, the balanced measure in Figure \ref{fig:path-graph} is basic vacuously, being unique; however, this can be shown more directly. If $\nu$ is any balanced measure with $\spt \nu$ equal to the endpoint set $\{ x_1, x_n \}$ and with $M_\nu = P_n$, then $(k-1) \+ \nu(x_1) + k \+ \nu(x_n)$ is constantly equal to the maximum of $T_\nu$ for $k = 1, ..., n$. In particular, taking $k = n$ gives $\nu(x_1)$ and taking $k = 1$ gives $\nu(x_n)$, so $\nu(x_1) = \nu(x_n)$ and all other vertices have measure $0$ by our hypothesis on the support. Since $\nu(P_n) = 1$, it follows that $\nu(x_1) = \tfrac{1}{2} = \nu(x_n)$, so $\nu = \mu$.

We also require a notion of compatibility between different probability measures.
\begin{defn} \label{defn:compatible}
    A family $\{ \mu_\alpha \}_{\alpha \in I}$ of probability measures on a graph is called \define{compatible} if
    \begin{equation*}
        \bigcup_{\alpha \in I} \spt \mu_\alpha \subseteq \bigcap_{\alpha \in I} M_{\mu_\alpha}.
    \end{equation*}
\end{defn}
This notion can be easily motivated as follows: suppose $\mu, \nu$ are two basic balanced measures. We can interpolate between these two balanced measures by considering $\lambda_t = (1-t) \+ \mu + t \+ \nu$ for $0 \leq t \leq 1$. Clearly $\spt \lambda_t = \spt \mu \cup \spt \nu$. Any vertex that maximizes both $T_\mu$ and $T_\nu$ necessarily maximizes $T_{\lambda_t}$ as well, so if the two measures are also compatible, then
\begin{equation*}
    \spt \lambda_t = \spt \mu \cup \spt \nu \subseteq M_\mu \cap M_\nu \subseteq M_{\lambda_t}.
\end{equation*}
Consequently, $\lambda_t$ is also balanced. Extending this idea to potentially larger families of measures leads to a notion of compatibility. This argument extends in a straight-forward way to linear combinations involving more than two balanced measures. Thus, a convex combination of compatible balanced measures is also a balanced measure. Our main result shows the converse result.

\begin{thm}[Decomposition of Balanced Measures] \label{thm:decomp}
    A measure on a graph is balanced if and only if it is a convex combination of compatible basic balanced measures.
\end{thm}

\begin{center}
\begin{figure}[h]
\begin{tikzpicture}
        \node at (0,0) {\includegraphics[width=3cm]{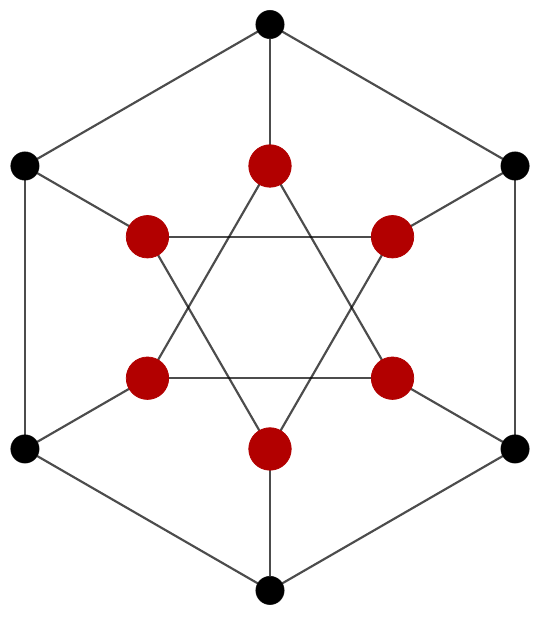}};
        \node at (6,0) {\includegraphics[width = 3cm]{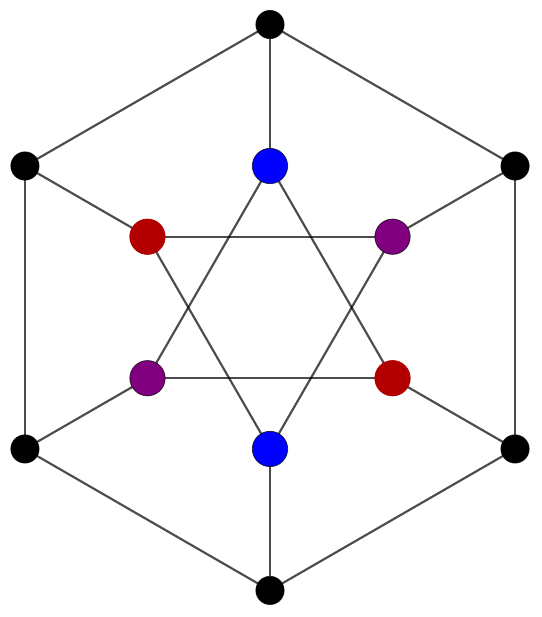}};
    \end{tikzpicture}
    \caption{Left: The support of a balanced measure on the D\"{u}rer graph. Right: A decomposition of this measure as a convex combination of three basic balanced measures, each of which is equidistributed on its support.}
    \label{fig:basic-nonminimal-petersen}
\end{figure}
\end{center}

Figure \ref{fig:basic-nonminimal-petersen} gives a simple example of a balanced measure and a corresponding decomposition into compatible basic balanced measures.

One way of interpreting the result is to think of the space of balanced measures as a finite union of polytopes in $\bbp(G)$, viewed as an $(n-1)$-simplex in the first orthant of $\R^n$ (see \S \ref{ss:notation}). The basic measures form the vertices of these polytopes, and the convex hulls of the compatible sets are the polytopes themselves.

In a more analytic vein, the basic balanced measures on $G$ are the balanced measures with a certain extremal property: namely, they are precisely those with the smallest possible cost-maximizing sets and the largest possible supports (cf.~Lemma \ref{lem:basic-iff}). In this setting, Theorem \ref{thm:decomp} states that every balanced measure on $G$ is a convex combination of these ``maximal" measures, and it prescribes exactly how to generate all balanced measures from this collection (viz., by convex combinations of \textit{compatible} subcollections). This generating set is moreover \textit{minimal} in the sense of Corollary \ref{cor:weak-minimality}, so Theorem \ref{thm:decomp} minimizes the search space in the problem of identifying all the balanced measures on $G$.

This naturally leads to the question of how many distinct balanced measures exist on a given graph. The example of the cycle graph $C_4$ shows that there can be an infinite number of balanced measures on a simple graph, but only two of these are basic. In fact, the number of basic balanced measures on a graph is at most exponential in the number of vertices, and this estimate is the best possible modulo the base of the exponent.

\begin{prop} \label{prop:bounds}
    The maximum number $b_n$ of basic balanced measures on a graph on $n \geq 2$ vertices satisfies
    \begin{equation*}
        2^{\lfloor \frac{n}{3} \rfloor}-1 \leq b_n \leq 2^{2n}-1.
    \end{equation*}
    Additionally, if $c_n$ is the maximum number of connected components in the space of balanced measures on a graph with $n \geq 2$ vertices, then
    \begin{equation*}
        2^{\lfloor \frac{n}{3} \rfloor}-1 \leq c_n \leq 2^{2n}-1.
    \end{equation*}
\end{prop}

\subsection{Motivation} \label{ss:motivation}

While the notion of balanced measure is new, it is directly related to a number of problems in very classical potential theory. Bj\"{o}rck \cite{bjorck1956distributions}, motivated by the classic problem from geometric measure theory of minimizing the \textit{$\alpha$-Riesz energy}
\begin{equation} \label{eq:riesz-energy}
    E_\alpha(\mu) := \int_{X \times X} \| x-y \|^{\alpha-n} \, d\mu(x) d\mu(y), \qquad 0 < \alpha \leq n,
\end{equation}
over all Borel probability measures supported on a given set $X \subseteq \R^n$, asked which Borel probability measures instead \textit{maximize} the energy
\begin{equation} \label{eq:lambda-energy}
    \medtilde{E}_\lambda(\mu) := \int_{X \times X} \| x-y \|^\lambda \, d\mu(x) d\mu(y), \qquad \lambda > 0.
\end{equation}
He showed that a maximizer always exists and is supported on the boundary $\partial X$. More abstractly, on a given metric space $(X,d)$, one might consider the problem of maximizing the value of
\begin{equation*}
    E(\mu) := \int_{X \times X} d(x,y) \, d\mu(x) d\mu(y)
\end{equation*}
over all Borel probability measures on $X$. The reader is referred to \cite{alexander1975generalized, alexander1977two, alexander1974extremal, carando2015energy, hinrichs2011note, wolf1997average, wolf1999averaging} for subsequent work in this direction. In the particular case that $X = G$ is a graph on $n$ vertices with distance matrix $D = \big( d(v_i, v_j) \big)_{i,j=1}^n$, then the energy functional defined by
\begin{equation} \label{eq:variational-prob}
    \mu \mapsto \langle \mu, D\mu \rangle
\end{equation}
is maximized over all $\mu \in \bbp(G)$ by a balanced measure (cf.~\cite{steinerberger2023sums} Proposition 1). In this sense, the study of balanced measures proceeds naturally from much older variational problems.

In a similar vein, one can relate this to a problem from spectral graph theory. The Courant--Fischer theorem states that the maximizers of the Rayleigh--Ritz quotient
\begin{equation*}
    Q_D(x) := \frac{\langle x, Dx \rangle}{\|x\|^2}, \qquad x \in \R^n \setminus \{0\},
\end{equation*}
are precisely the eigenvectors corresponding to the largest eigenvalue of $D$. Replacing the Euclidean norm $\| \cdot \|$ with the $1$-norm $\| \cdot \|_1$ and applying a change of parameters $x/\|x\|_1 \mapsto \mu$, one sees that the problem in Equation \eqref{eq:variational-prob} is precisely the $1$-norm analogue of the above, but restricted to the first orthant.

Another motivation is a 1964 result of Gross \cite{gross1964rendezvous}: for every compact, connected metric space $(X,d)$, there exists a unique number $r(X,d) \in \R$ such that, for every finite set $\{ x_1, \dots, x_N \} \subset X$, there is another point $x \in X$ satisfying
\begin{equation*}
    \frac{1}{N} \sum_{i=1}^N d(x,x_i) = r(X,d).
\end{equation*}
The survey by Cleary, Morris, \& Yost \cite{cleary1986numerical} gives a wonderful introduction to the subject. Note that no such results can be true on graphs, but see Thomassen \cite{thomassen2000rendezvous} for a substitute result. (See also \cite{steinerberger2023curvature} for a connection to the von Neumann minimax theorem). The result of Gross can be extended to more general probability measures, a consequence of which is the following: if there exists $\mu \in \bbp(X)$ such that the assignment
\begin{equation*}
    x \mapsto \int_X d(x,y) \, d\mu(y)
\end{equation*}
is constant, then that constant is precisely $r(X,d)$. This particular condition is the Euler--Lagrange equation related to the problem of maximizing $E(\mu)$. This is one way allowing for an explicit computation of $r(X,d)$ which, in general, is very difficult. We also note a philosophical connection to earlier work in \cite{larcher2000approximation} that proceeds along similar lines in quasihypermetric spaces, which itself is a direct outgrowth of the aforementioned work of Bj\"{o}rck \cite{bjorck1956distributions}. The structural theorems in \cite{larcher2000approximation} on the limits of the greedy sequences of measures are not dissimilar to our own Theorem \ref{thm:decomp}.

\subsection{Notation} \label{ss:notation}

Throughout this paper, a graph is always assumed to be connected, simple, and finite with unit edge weights. We frequently identify the space $\bbp(G)$ of all probability measures on an $n$-vertex graph $G$ with an $(n-1)$-simplex in $\R^n$:
\begin{equation*}
    \bbp(G) \cong \left\{ (x_1, ..., x_n) \in \R^n \!: x_i \geq 0 \ \forall \+ i \text{ and } \sum_{i=1}^n x_i = 1 \right\},
\end{equation*}
where if $V = \{ v_1, ..., v_n \}$ for some fixed enumeration of the vertex set, then $\mu \cong (\mu(v_1), ..., \mu(v_n))$. This viewpoint affords another way of phrasing the definition given in Equation \eqref{eq:balanced-defn}. Denote by $D$ the distance matrix of $G$. Then $D\mu = \big( T_\mu(v_1), ..., T_\mu(v_n) \big)$, so $\mu$ is balanced precisely when
\begin{equation*}
    v_i \in \spt \mu \quad \text{implies} \quad (D\mu)_i = \| D\mu \|_\infty,
\end{equation*}
where $(D\mu)_i$ is the $i$th coordinate of $D\mu$.

The two subfamilies of $\bbp(G)$ that we study are the family $\bbb(G) \subseteq \bbp(G)$ of \textit{balanced} measures and the family $\mathsf{B}(G) \subseteq \bbb(G)$ of \textit{basic balanced} measures. Lastly, we recall from \S \ref{ss:balanced} that $\spt \mu$ denotes the \textit{support} of a measure on a graph, $T_\mu$ its \textit{transport cost function}, and $M_\mu$ its \textit{cost-maximizing set}.

\subsection{Outline of the Paper}

Our main purpose is to prove Theorem \ref{thm:decomp} and Proposition \ref{prop:bounds}; this is done in \S \ref{s:basic-measures} and \S \ref{s:joins}, respectively. Section \ref{s:basic-measures} begins with several results on compatible collections of balanced measures, after which the key ``one-sided extrapolation" Lemma \ref{lem:one-sided} is proved. From there, the proof of Theorem \ref{thm:decomp} comes down to a judicious choice of notation. Rephrasing the result in terms of a relation on $\mathcal{P}(G) \times \mathcal{P}(G)$ that compares the supports and cost-maximizing sets of pairs of measures, we can then readily apply the one-sided extrapolation lemma and establish the main theorem. The section concludes with a discussion of the extent to which the basic balanced measures form a minimal collection from which the other balanced measures on $G$ can be generated through convex combinations.

Subsequently in \S \ref{s:joins}, we study the interplay between balanced measures and the graph join operation. With a few technical results in hand, the lower bound in Proposition \ref{prop:bounds} readily falls out of a graph join construction, whereas the upper bound is a simple consequence of Theorem \ref{thm:decomp}. Finally, \S \ref{s:examples} demonstrates the richness of the structure of the family of balanced measures through several examples. In particular, we show in Proposition \ref{prop:subgraph} that there is a natural graph structure on the family $\mathsf{B}(G)$ of basic balanced measures, and that every graph can be realized as a subgraph of some $\mathsf{B}(G)$ in a natural way. This section is essentially self-contained and may be read independently of \S \ref{s:basic-measures}--\ref{s:joins}.

\section{Minimal Generating Sets of the Space of Balanced Measures} \label{s:basic-measures}

\subsection{Overview} In this section, we prove Theorem \ref{thm:decomp} and explore some of its immediate consequences. We begin by showing the ``easy direction," i.e., that convex combinations of compatible balanced measures are balanced. This is followed by a sort of converse result that is utilized near the end of the section.  

We then discuss the central Lemma \ref{lem:one-sided} on extrapolating between a pair of balanced measures. This result is then recontextualized by a partial ordering we place on the set of pairs $(\spt \mu, M_{\mu})$ associated with balanced measures $\mu$. With this relation in hand, we show that the basic measures are not only precisely the balanced measures that yield maximal pairs with respect to this partial ordering, but also that, whenever $\mu$ is in the set $\bbb(G)$ of balanced measures, the pair $(\spt \mu, M_{\mu})$ is upper bounded by $(\spt \nu, M_{\nu})$ for some \textit{basic} balanced $\nu \in \mathsf{B}(G)$. These final results are the last pieces needed to prove the more difficult direction of Theorem \ref{thm:decomp}. The section with a discussion of the weak sense in which the set of basic measures ``minimally" generates the space of balanced measures.

\subsection{Proof of Theorem \ref{thm:decomp}}

The first lemma in the proof of Theorem \ref{thm:decomp} is simply the ``if" direction.

\vspace*{-2pt}

\begin{lem} \label{lem:convex-of-compatible}
    Let $\{ \mu_i \}_{i=1}^N$ be a compatible set of balanced measures on a graph $G = (V,E)$ and $\mu = \sum_{i=1}^N a_i \+ \mu_i$ a convex combination thereof. Then $\mu$ is balanced.
\end{lem}

\vspace*{-2pt}

\pf Since the entries of each $\mu_i$ sum to $1$, the entries of each $a_i \+ \mu_i$ sum to $a_i$. In turn, $\sum_{i=1}^N a_i = 1$ so the entries of $\mu$ sum to $1$ as well. As each $a_i \geq 0$ and each $\mu_i$ takes nonnegative values, $\mu$ is a probability measure on $G$.

We now show that $\mu$ is balanced. Let $S = \bigcup_{i=1}^N \spt \mu_i$ and $M = \bigcap_{i=1}^N M_{\mu_i}$, and recall that $S \subseteq M$ because $\{ \mu_i \}_{i=1}^N$ is compatible. Certainly $\mu$ is supported on $S$, as every $\mu_i$ is identically $0$ outside of $S$. Hence, $\spt \mu \subseteq M$, and for all $v \in \spt \mu$ and $w \in V$ we have
\begin{equation*}
    T_{\mu}(v) = \sum_{i=1}^N a_i T_{\mu_i}(v) \geq \sum_{i=1}^N a_i T_{\mu_i}(w) = T_{\mu}(w)
\end{equation*}
where the inequality holds because $v$ is in the cost-maximizing set of every $\mu_i$. This shows that $v \in M_{\mu}$, so $\spt \mu \subseteq M_{\mu}$; that is, $\mu$ is balanced. \textqed

Before building up results to show the other direction of the decomposition theorem\textemdash that \textit{all} balanced measures on a graph can be written as such a convex combination\textemdash we pause to shed light on the role of compatibility. In fact, compatibility is a \textit{necessary} condition for more than one (nontrivial) convex combination of two balanced measures to be balanced. In addition to finding use in the proof of Theorem \ref{thm:decomp}, this result is instructive and interesting in its own right.

\vspace*{-2pt}

\begin{lem} \label{lem:if-balanced}
    Let $\mu,\nu \in \bbb(G)$. If the measure $\lambda_t := (1-t) \+ \mu + t \+ \nu$ is balanced for two distinct $t \in (0,1)$, then $\mu$ and $\nu$ are compatible.
\end{lem}

\vspace*{-2pt}

\pf Write the transport cost function $T_{\lambda_t}$ as
\begin{equation*}
    T_{\lambda_t} = (1-t) \+ T_\mu + t \+ T_\nu.
\end{equation*}
Suppose $\lambda_t$ is balanced for distinct values $t = t_1, t_2 \in (0,1)$ and let $u \in \spt \mu \subseteq M_\mu$ and $v \in \spt \nu \subseteq M_\nu$. Then $u,v \in \spt \lambda_t$ as well, so for $t = t_1, t_2$ we have
\begin{align*}
    \max_{w \in V} &\ T_{\lambda_t}(w) = T_{\lambda_t}(u) \\
    &= (1-t) T_{\mu}(u) + t \+ T_{\nu}(u) \\[0.1cm]
    &= (1-t) \max_{w \in V} T_{\mu}(w) + t T_{\nu}(u)
\end{align*}
and, similarly,
\begin{equation*}
    \max_{w \in V} T_{\lambda_t}(w) = (1-t) T_{\mu}(v) + t \max_{w \in V} T_{\nu}(w).
\end{equation*}
It follows that
\begin{equation*}
    (1-t) \Big( \max_{w \in V} T_{\mu}(w) - T_{\mu}(v) \Big) = t \Big( \max_{w \in V} T_{\nu}(w) - T_{\nu}(u) \Big)
\end{equation*}
for two distinct $t \in (0,1)$, which forces
\begin{equation*}
    \max_{w \in V} T_{\mu}(w) - T_{\mu}(v) = 0 \quad \text{and} \quad \max_{w \in V} T_{\nu}(w) - T_{\nu}(u) = 0.
\end{equation*}
Thus, $v \in M_{\mu}$ and $u \in M_{\nu}$, so $\spt \mu \cup \spt \nu \subseteq M_{\mu} \cap M_{\nu}$; i.e., $\mu$ and $\nu$ are compatible. \textqed

We combine this with Lemma \ref{lem:convex-of-compatible} into a single result.

\begin{cor} \label{cor:convex-combo-iff}
    If $\mu,\nu \in \bbb(G)$, then the measure $\lambda_t := (1-t) \+ \mu + t \+ \nu$ is balanced for all $t \in (0,1)$ if and only if $\mu$ and $\nu$ are compatible.
\end{cor}

Likewise, Lemma \ref{lem:if-balanced} combines with Lemma \ref{lem:convex-of-compatible} to give a characterization of the finite sets in $\bbb(G)$ whose convex hulls are contained in $\bbb(G)$.

\begin{cor} \label{cor:compatible-iff}
    Let $\{ \mu_i \}_{i=1}^N$ be a finite set of balanced measures on a graph. Then all convex combinations of measures in $\{ \mu_i \}_{i=1}^N$ are balanced if and only if the collection $\{ \mu_i \}_{i=1}^N$ is compatible.
\end{cor}

\pf If $\{ \mu_i \}_{i=1}^N$ is compatible, then all convex combinations are balanced by Lemma \ref{lem:convex-of-compatible}. Conversely, if $\mathcal{M}$ is not compatible, then some vertex $v \in \bigcup_{i=1}^N \spt \mu_i$ does not belong to $\bigcap_{i=1}^N M_{\mu_i}$, the intersection of the cost-maximizing sets. Thus, there are (at least) two distinct $\mu_j,\mu_k$ such that $v \in \spt \mu_j \cup \spt \mu_k$ but $v \not\in M_{\mu_j} \cap M_{\mu_k}$. By Corollary \ref{cor:convex-combo-iff}, there exist $a_j = (1-t), a_k = t \in (0,1)$ for which
\begin{equation*}
    \sum_{i=1}^N a_i \+ \mu_i = (1-t) \+ \mu_j + t \+ \mu_k
\end{equation*}
is not balanced, where $a_i = 0$ for all $i \neq j,k$. In particular, not all convex combinations of measures in $\{ \mu_i \}_{i=1}^N$ are balanced. \textqed

This concludes the ``if" direction of Theorem \ref{thm:decomp}. We undertake the more challenging ``only if" direction beginning with the following technical lemma.

\begin{lem} \label{lem:interp-sup-max}
    Let $G = (V,E)$ be a graph and let $\mu,\nu$ be a pair of compatible balanced measures. Then, for all $t \in (0,1)$ the measure $\lambda_t := (1-t) \+ \mu + t \+ \nu$ is supported on $\spt \mu \cup \spt \nu$ and maximized on $M_\mu \cap M_\nu$: that is,
    \begin{equation*}
        \spt \lambda_t = \spt \mu \cup \spt \nu \quad \text{and} \quad M_{\lambda_t} = M_\mu \cap M_\nu \qquad \forall \+ t \in (0,1).
    \end{equation*}
\end{lem}

\pf Let $t \in (0,1)$. Then both $t > 0$ and $1-t > 0$, and if $v \in \spt \mu \cup \spt \nu$, then at least one of $\mu(v)$, $\nu(v)$ is positive. Hence,
\begin{equation*}
    \lambda_t(v) = (1-t) \+ \mu(v) + t \+ \nu(v) > 0.
\end{equation*}
On the other hand, if $v \not\in \spt \mu \cup \spt \nu$, then both $\mu(v) = 0$ and $\nu(v) = 0$, whence $\lambda_t(v) = 0$.

Now, since $T_\mu$ is constant on $M_\mu$ and $T_\nu$ is constant on $M_\nu$, the cost function $T_{\lambda_t} = (1-t) T_{\mu} + t T_{\nu}$ is constant on $M_\mu \cap M_\nu$. Furthermore, if $v \in M_\mu \cap M_\nu$ and $w \not\in M_\nu$, then $T_\mu(v) \geq T_\mu(w)$ and $T_\nu(v) > T_\nu(w)$, so
\begin{equation*}
    T_{\lambda_t}(v) = (1-t) T_\mu(v) + t T_\nu(v) > (1-t)T_\mu(w) + tT_\nu(w) = T_{\lambda_t}(w).
\end{equation*}
By symmetry in $\mu$ and $\nu$, the same conclusion holds if instead $w \not\in M_\mu$, meaning $T_{\lambda_t}$ is not maximized at any point outside of $M_{\mu} \cap M_{\nu}$. As $T_{\lambda_t}$ is constant on $M_{\mu} \cap M_{\nu}$, we have that $T_{\lambda_t}(v) \geq T_{\lambda_t}(w)$ for all $v \in M_{\mu} \cap M_{\nu}$ and $w \in V$, with equality if and only if $w \in M_{\mu} \cap M_{\nu}$. This proves that $M_{\mu} \cap M_{\nu}$ is the cost-maximizing set of $\lambda_t$. \textqed

This springboards us into the following key lemma, which analyzes the \textit{linear} combinations of a pair of balanced measures, yielding conclusions about their supports and cost-maximizing sets in the case that these measures are balanced. It is through extrapolation rather than interpolation\textemdash that is, by taking linear combinations with coefficients outside the interval $[0,1]$\textemdash that we identify the boundary of $\bbb(G)$. Doing so is a crucial first step toward identifying the basic balanced measures in $\bbb(G)$.

\begin{lem}[One-Sided Extrapolation Lemma] \label{lem:one-sided}
    Let $\mu$ and $\nu$ be distinct balanced measures on a graph $G = (V,E)$ with $\spt \mu \subseteq \spt \nu$ and $M_\mu \supseteq M_\nu$. For each $t \in \R$, let
    \begin{equation*}
        \lambda_t := (1-t) \+ \mu + t \+ \nu.
    \end{equation*}
    Then the set
    \begin{equation*}
        \{ t \in \R \!: \lambda_t \text{ is a balanced measure} \+ \}
    \end{equation*}
    is a compact interval $[L,R]$ with $L \leq 0$ and $R > 1$. Furthermore, we have the following inclusions
    \begin{equation} \label{eq:lambda-R}
        \spt \lambda_R \subseteq \spt \nu \quad \text{and} \quad M_{\lambda_R} \supseteq M_\nu,
    \end{equation}
    where at least one of the inclusions is strict, and we have the equalities
    \begin{equation} \label{eq:lambda-t-nu}
        \spt \lambda_t = \spt \nu \quad \text{and} \quad M_{\lambda_t} = M_\nu \qquad  \forall \+ t \in (L,R).
    \end{equation}
\end{lem}

\pf \textsc{Step 1.} We show that the set of $t \in \R$ for which $\lambda_t \in \bbp(G)$ is a closed, bounded interval. For each $v \in V$, let $I_v := \{ t \in \R \!: \lambda_t(v) \geq 0 \}$. Since $t \mapsto (1-t) \+ \mu(v) + t \+ \nu(v)$ is continuous and monotonic, it must be that $I_v$ is a closed, unbounded interval of the form
\begin{equation} \label{eq:closed-intervals}
    I_v = \left\{ \begin{array}{cl}
        [a_v,\infty) & \text{if } \ \mu(v) < \nu(v) \\
        \R & \text{if } \ \mu(v) = \nu(v) \\
        (-\infty,b_v] & \text{if } \ \mu(v) > \nu(v).
    \end{array} \right.
\end{equation}
Notice that
\begin{equation} \label{eq:measure-iff-nonneg}
\begin{aligned}
    \bigcap_{v \in V} I_v &= \{ t \in \R \!: \lambda_t(v) \geq 0 \ \forall \+ v \in V \} \\
    &= \{ t \in \R \!: \lambda_t \text{ is a measure} \+ \} \\
    &= \{ t \in \R \!: \lambda_t \text{ is a probability measure} \+ \},
\end{aligned}
\end{equation}
as the only requirement for $\lambda_t$ to be a measure is that it have nonnegative entries, and the condition $\sum_{v \in V} \lambda_t(v) = 1$ is always satisfied because
\begin{equation*}
    \sum_{v \in V} \big( (1-t) \+ \mu(v) + t \+ \nu(v) \big) = (1-t) \sum_{v \in V} \mu(v) + t \sum_{v \in V} \nu(v) = (1-t) + t = 1.
\end{equation*}
Now we claim that $\bigcap_{v \in V} I_v$ is a compact interval $[L_0,R_0]$. Since it is a finite intersection of closed intervals, it is enough to show that at least one of the $I_v$ is left-bounded and at least one is right-bounded. Recall that $\mu$ and $\nu$ are distinct. Since their respective values sum to $1$, neither can be greater than or equal to the other at all points of $V$, so there exist vertices $v_0, v_1 \in V$ such that $\mu(v_0) < \nu(v_0)$ and $\mu(v_1) > \nu(v_1)$. By Equation \eqref{eq:closed-intervals}, $I_{v_0} = \big[ a_{v_0}, \infty \big)$ and $I_{v_1} = \big( -\infty, b_{v_1} \big]$, so we are done. In fact, $0,1 \in \bigcap_{v \in V} I_v$ by Equation \eqref{eq:measure-iff-nonneg}, so this intersection contains the entire interval $[0,1]$. Notice also that $\bigcap_{v \in V} I_v = \bigcap_{v \in \spt \nu} I_v$, because $\mu(w) = \nu(w) = 0$ for $w \not\in \spt \nu \supseteq \spt \mu$, in which case Equation \eqref{eq:closed-intervals} gives that $I_w = \R$.

\noindent \textsc{Step 2.} Having shown that the parameters $t$ for which $\lambda_t$ is a probability measure form a closed, bounded interval $[L_0,R_0] \supseteq [0,1]$, we show that the same is true of the subset consisting of those $t$ for which $\lambda_t$ is balanced. For each $w \in V$, let $J_w$ be the set of parameters $t$ such that the transport cost of $\lambda_t$ at $w$ is no greater than the transport cost on $M_\nu$:
\begin{equation*}
    J_w := \big\{ t \in \R \!: T_{\lambda_t}(v) \geq T_{\lambda_t}(w) \ \forall \+ v \in M_\nu \big\}.
\end{equation*}
(Since $\nu$ is constant on $M_\nu$, we in fact have $t \in J_w$ if and only if $T_{\lambda_t}(v) \geq T_{\lambda_t}(w)$ for \textit{some} $v \in M_\nu$.) Notice that the maximizing set $M_{\lambda_t}$ contains $M_\nu$ if and only if $t \in J_w$ for all $w \in V \setminus M_\nu$. Because $\spt \mu \subseteq \spt \nu$ by hypothesis and $\spt \nu \subseteq M_\nu$ by the definition of a balanced measure, we have that $\spt \lambda_t \subseteq M_\nu$ for all $t$, so
\begin{equation} \label{eq:balanced-params}
\begin{aligned}
    \bigcap_{v \in \spt \nu} I_v \cap \bigcap_{w \in V \setminus M_\nu} J_w &= \big\{ t \in \R \!: \lambda_t \text{ is a probability measure and } M_{\lambda_t} \subseteq M_\nu \big\} \\
    &= \{ t \in \R \!: \lambda_t \text{ is a balanced measure} \+ \}.
\end{aligned}
\end{equation}
By reasoning identical to that in Step 1, the intervals $J_w$ take the form
\begin{equation*}
    J_w = \left\{ \begin{array}{cl}
        [\alpha_w,\infty) & \text{if } \ T_{\mu}(w) < T_{\nu}(w) \\
        \R & \text{if } \ T_{\mu}(w) = T_{\nu}(w) \\
        (-\infty,\beta_w] & \text{if } \ T_{\mu}(w) > T_{\nu}(w),
    \end{array} \right.
\end{equation*}
so the set in Equation \eqref{eq:balanced-params} is a closed, bounded interval that we call $[L,R]$. As before, we observe that $0,1 \in [L,R]$, so $[L,R] \supseteq [0,1]$.

\noindent \textsc{Step 3.} We prove that in fact $R > 1$. In Step 2, we have already seen that $\spt \lambda_t \subseteq M_\nu$ and $M_{\lambda_t} \supseteq M_\nu$ for all $t$ such that $\lambda_t \in \bbb(G)$, i.e., for all $\lambda_t \in [L,R]$. In particular, \eqref{eq:lambda-R} holds, and we claim that at least one of these relations is strict; the claim that $R > 1$ will readily follow.

Since $R$ is an endpoint of a finite intersection of closed intervals, it must be an endpoint of one of the intervals themselves. The intervals $I_v$ and $J_w$ are the positive sets and negative sets of the affine functions
\begin{gather*}
    t \mapsto (1-t) \+ \mu(v) + t \+ \nu(v) \qquad \text{and} \\
    t \mapsto (1-t) \left( \max_{u \in V} T_{\mu}(u) - T_{\mu}(w) \right) + t \left( \max_{u \in V} T_{\nu}(u) - T_{\nu}(w) \right),
\end{gather*}
so by continuity, their endpoints are points at which those functions vanish. In particular, if $R$ is an endpoint of $I_v$ for some $v \in \spt \nu$, then $\lambda_R(v) = 0$, so $\lambda_R$ vanishes at some point of $\spt \nu$; that is, $\spt \lambda_R \subsetneq \spt \nu$. Likewise, if instead $R$ is an endpoint of $J_w$ for some $w \in V \setminus M_{\nu}$, then $T_{\lambda_R}$ attains its maximum at $w$. We saw already that $T_{\lambda_t}$ is also maximized on $M_{\nu}$, and since $w \not\in M_{\nu}$, we have $M_{\lambda_R} \supsetneq M_{\nu}$. In both cases, we have either $\spt \lambda_R \neq \spt \nu$ or $M_{\lambda_R} \neq M_\nu$, so $\lambda_R \neq \nu = \lambda_1$. Recalling that $[0,1] \subseteq [L,R]$, we conclude that $R > 1$.

\noindent \textsc{Step 4.} It remains to show that for $t \in (L,R)$, $\spt \lambda_t = \spt \nu$ and $M_{\lambda_t} = M_{\nu}$. It is clear that, for each $t \in (L,R)$, the measure $\lambda_t$\textemdash being a linear combination of $\mu$ and $\nu$\textemdash is also a convex combination of $\lambda_L$ and $\lambda_R$. By Lemma \ref{lem:if-balanced}, since all such $\lambda_t$ are balanced, $\lambda_L$ and $\lambda_R$ must be compatible. It then follows from Lemma \ref{lem:interp-sup-max} that the $\lambda_t$ have the same support and cost-maximizing sets for all $t \in (L,R)$. Since $1 \in (L,R)$ and $\lambda_1 = \nu$, the equality of these $\spt \lambda_t$ and $M_{\lambda_t}$ implies Equation \eqref{eq:lambda-t-nu}, completing the proof of the lemma. \textqed

This motivates a poset structure in terms of which Lemma \ref{lem:one-sided} admits a more succinct restatement.

\begin{defn} \label{defn:poset}
    Let $G = (V,E)$ be a graph and let $\mathcal{P}(V)$ be the power set of $V$. Define the following poset relation on $\mathcal{P}(V) \times \mathcal{P}(V)$: let
    \begin{equation*}
        (S,M) \leq (S',M') \quad \text{if and only if} \quad S \subseteq S' \text{ and } M \supseteq M'.
    \end{equation*}
\end{defn}

\begin{cor} \label{cor:one-sided-rephrase}
    Let $G = (V,E)$ be a graph and let $\mu,\nu$ be two measures such that $(\spt \mu, M_{\mu}) \geq (\spt \nu, M_{\nu})$. Then there exists a measure $\lambda$ such that $(\spt \lambda, M_{\lambda}) > (\spt \nu, M_{\nu})$ and $\nu$ is a convex combination of $\lambda$ and $\mu$. There also exists an infinite family $S \subseteq \bbb(G)$ such that, for each $\rho \in S$, $(\spt \rho, M_{\rho}) = (\spt \nu, M_{\nu})$.
\end{cor}

\pf The hypothesis of the corollary is identical to the hypothesis of Lemma \ref{lem:one-sided}. We take the $\lambda$ attested to above to be the $\lambda_R$ of Lemma \ref{lem:one-sided}, and $S$ to be the set $\{ \lambda_t : t \in (L,R) \}$. Then, except for the statement that $\nu$ is a convex combination of $\lambda_R$ and $\mu$, all the asserted conclusions follow immediately from Lemma \ref{lem:one-sided} and the definition of the partial ordering. To prove that $\nu$ is such a convex combination, we simply use the definition $\lambda_R = (1-R) \+ \mu + R \+ \nu$, from which we see that
\begin{equation*}
    \frac{1}{R} \+ \lambda_R + \frac{R-1}{R} \+ \mu = \frac{1}{R} \big( (1-R) \+ \mu + R \+ \nu \big) + \frac{R-1}{R} \mu = \nu.
\end{equation*}
We have that $1/R + (R-1)/R = 1$ and $R > 1$ by Lemma \ref{lem:one-sided}, so $1/R, (R-1)/R > 0$. This presents $\nu$ as a convex combination of $\lambda_R$ and $\mu$. \textqed

This leads to an equivalent characterization of basic measures.

\begin{lem} \label{lem:basic-iff}
    Let $G = (V,E)$ be a graph and $\nu \in \bbb(G)$. Then $\nu$ is basic if and only if $(\spt \nu, M_\nu)$ is maximal in $\{ (\spt \rho, M_{\rho}) : \rho \text{ is a balanced measure} \+ \}$ with respect to the partial ordering of Definition \ref{defn:poset}.
\end{lem}

\pf We prove that the negations of these conditions are equivalent. Assume that $\nu$ is not basic. Then there exists $\mu \in \bbp(G) \setminus \{ \nu \}$ such that $\spt \mu = \spt \nu$ and $M_{\mu} = M_{\nu}$. This pair of measures satisfies the hypotheses of Corollary \ref{cor:one-sided-rephrase}, so there exists $\lambda \in \bbb(G)$ such that $(\spt \lambda, M_{\lambda}) > (\spt \nu, M_{\nu})$ by the corollary. Thus, $(\spt \nu, M_{\nu})$ is not maximal.

Assume conversely that $(\spt \nu, M_{\nu})$ is not maximal. Then there exists $\mu \in \bbp(G)$ such that $(\spt \mu, M_{\mu}) > (\spt \nu, M_{\nu})$. Then, by Corollary \ref{cor:one-sided-rephrase}, there exists an infinite family $S \subseteq \bbb(G)$ such that $(\spt \rho, M_{\rho}) = (\spt \nu, M_{\nu})$ for all $\rho \in S$. The existence of any such $\rho \neq \nu$ implies, by definition, that $\nu$ is not basic. \textqed

In fact, we can further conclude that each element in this poset is upper bounded by the pair given by some basic measure.

\begin{lem} \label{lem:basic-gt-balanced}
    Let $G = (V,E)$ be a graph and let $\rho$ be a balanced measure on it. Then there exists some basic measure $\mu$ such that $(\spt \mu, M_\mu) \geq (\spt \rho,M_\rho)$.
\end{lem}

\pf Let $S = \{ (\spt \sigma, M_\sigma) : \sigma \text{ is a balanced measure} \+ \}$. If $(\spt \rho, M_{\rho})$ is not maximal in $S$, there exists a balanced measure $\rho_1$ such that $(\spt \rho, M_\rho) < (\spt \rho_1, M_{\rho_1})$. In turn, if $(\spt \rho_1, M_{\rho_n})$ is not maximal in $S$, there exists a balanced measure $\rho_2$ such that $(\spt \rho_1, M_{\rho_1}) < (\spt \rho_2,M_{\rho_2})$. We may continue this process, forming a chain of elements in the poset:
\begin{equation*}
    (\spt \rho, M_\rho) < (\spt \rho_1, M_{\rho_1}) < (\spt \rho_2,M_{\rho_2}) < \cdots.
\end{equation*}
This process must eventually terminate because $S$ is finite, so we have a finite chain
\begin{equation*}
    (\spt \rho, M_\rho) < (\spt \rho_1, M_{\rho_1}) < (\spt \rho_2, M_{\rho_2}) < \cdots < (\spt \rho_k, M_{\rho_k})
\end{equation*}
where $(\spt \rho_k, M_{\rho_k})$ is maximal in $S$. Then $\rho_k$ is basic by Lemma \ref{lem:basic-iff}, which completes the proof. \textqed

We are at last ready to complete the proof of Theorem \ref{thm:decomp}.

\noindent \define{Proof of Theorem \ref{thm:decomp}.} We have shown in Lemma \ref{lem:if-balanced} that a convex combination of compatible basic balanced measures is balanced, so it remains to prove the following statement:
\begin{quote}
    \textit{If $\rho \in \bbb(G)$ is a balanced measure on $G$, then $\rho$ is a convex combination of a compatible set of basic balanced measures on $G$.}
\end{quote}
To accomplish this, we construct a finite sequence of balanced measures on $G$, starting with $\rho_1 = \rho$. By Lemma \ref{lem:basic-gt-balanced}, there exists some basic balanced measure $\mu_1 \in \mathsf{B}(G)$ such that $\big( \spt \mu_1, M_{\mu_1} \big) \geq \big( \spt \rho_1, M_{\rho_1} \big)$. If $\rho_1$ is not basic, then $\mu_1 \neq \rho_1$ and, by Corollary \ref{cor:one-sided-rephrase}, there exists a balanced measure $\rho_2$ such that $\big( \spt \rho_1, M_{\rho_1} \big) < \big( \spt \rho_2, M_{\rho_2} \big)$ and $\rho_1$ is a convex combination of $\rho_2$ and $\mu_1$.

If $\rho_2$ is not basic, there exists some basic $\mu_2$ such that $\big( \spt \mu_2, M_{\mu_2} \big) \geq \big( \spt \rho_2,$ $M_{\rho_2} \big)$ and again by Corollary \ref{cor:one-sided-rephrase}, there exists a balanced measure $\rho_3$ such that $\big( \spt \rho_2, M_{\rho_2} \big) < \big( \spt \rho_3, M_{\rho_3} \big)$ and $\rho_2$ is a convex combination of $\rho_3$ and $\mu_2$.

Just as in the previous lemma, we form an eventually terminating chain
\begin{equation*}
    (\spt \rho, M_\rho) < \big( \spt \rho_1, M_{\rho_1} \big) < \big( \spt \rho_2, M_{\rho_2} \big) < \cdots < \big( \spt \rho_k, M_{\rho_k} \big),
\end{equation*}
where $\rho_k$ is basic (or else the chain would not have terminated) and, for each $j < k$, $\rho_t$ is a convex combination of $\rho_{j+1}$ and a balanced measure $\mu_j$ such that $\big( \spt \rho_j, M_{\rho_j} \big) < \big( \spt \mu_j, M_{\mu_j} \big)$.

We then prove by finite descent on the indices that, for each $j = k, ..., 1$, the measure $\rho_j$ is a convex combination measures in the set $S = \{ \mu_1, \mu_2, ..., \mu_{k-1}, \rho_k \}$. In the base case, $\rho_k$ is such a combination because it is in $S$. In the inductive step, we assume that $\rho_{j+1}$ is such a combination and must prove that $\rho_j$ is. But $\rho_j$ is a convex combination of $\rho_{j+1}$ and $\mu_j$, so this is also clear by the inductive hypothesis.

Thus in particular, $\rho_1$ is a convex combination of the set $S = \{ \mu_1, \mu_2, ..., \mu_{k-1}, \rho_k \}$ of basic measures. We also see that $\big( \spt \rho_1, M_{\rho_1} \big) \leq \big( \spt \rho_k, M_{\rho_k} \big)$, and for each $j$,
\begin{equation*}
    (\spt \rho_1,M_{\rho_1}) \leq (\spt \rho_j, M_{\rho_j} \big) < \big( \spt \mu_j, M_{\mu_j} \big).
\end{equation*}
This shows that, for each $\sigma \in S$, $\spt \sigma \subseteq \spt \rho_1 \subseteq M_{\rho_1} \subseteq M_\sigma$ holds; hence,
\begin{equation*}
    \bigcup_{\sigma \in S} \spt \sigma \subseteq \spt \rho_1 \subseteq M_{\rho_1} \subseteq \bigcap_{\sigma \in S} M_\sigma.
\end{equation*}
Therefore, $S$ is a compatible set of basic balanced measures and $\rho = \rho_1$ is a convex combination of thereof, as we sought to show. \textqed

As indicated in Remark \ref{rmk:not-minimal} below, this set of generators is typically not minimal up to generating all balanced measures through convex combinations. However, it \textit{is} a minimal generating set in a weaker, somewhat more technical sense. To demonstrate this minimality, we introduce the following lemma.

\begin{lem} \label{lem:convex-combo-leq}
    Let $G = (V,E)$ be a graph and $\{ \mu_i \}_{i=1}^N \subseteq \bbb(G)$ a compatible set of basic balanced measures. In addition, let $\mu = \sum_{i=1}^N a_i \+ \mu_i$ be a balanced measure such that each $a_i$ is strictly positive. Then $(\spt \mu, M_\mu) \leq \big( \spt \mu_i, M_{\mu_i} \big)$ for all $i \in \{ 1, ..., N \}$.
\end{lem}

\pf Let $i \in \{ 1, ..., N \}$ be arbitrary and let $v \in \spt \mu_i$. Since the coefficients $a_j$ are strictly positive, we have $a_i \+ \mu_i(v) > 0$ and, for all $j \neq i$, $a_j \mu_j(v) \geq 0$; hence, $\mu(v) > 0$. That is, $v \in \spt \mu$, so $\spt \mu_i \subseteq \spt \mu$.

Now let $v \not\in M_{\mu_i}$. Let $M = \bigcap_{j=1}^N M_{\mu_j}$ and note that $M \supseteq \bigcup_{i=1}^N \spt \mu_i \neq \varnothing$, so there exists some $w \in M$. In particular $w \in M_{\mu_i}$, so $T_{\mu_i}(w) > T_{\mu_i}(v)$ and thus $a_i T_{\mu_i}(w) > a_i T_{\mu_i}(v)$ by our hypothesis that $a_i > 0$. At the same time, for $j \neq i$, $w \in M_{\mu_j}$ so $T_{\mu_j}(w) \geq T_{\mu_j}(v)$ and thus
\begin{equation*}
    \sum_{\substack{1 \leq j \leq N \\ j \neq i}} a_j T_{\mu_j}(w) \geq \sum_{\substack{1 \leq j \leq N \\ j \neq i}} a_j T_{\mu_j}(v).
\end{equation*}
Adding the inequality $a_i T_{\mu_i}(w) > a_i T_{\mu_i}(v)$ here yields
\begin{equation*}
    \mu(w) = \sum_{j=1}^N a_j T_{\mu_j}(w) > \sum_{j=1}^N a_j T_{\mu_j}(v) = \mu(v).
\end{equation*}
Therefore, $v \not\in M_\mu$, so we conclude that $M_\mu \subset M_{\mu_i}$. We have shown that $\spt \mu_i \subseteq \spt \mu$, and these two claims are enough to conclude the statement of the lemma by the definition of the poset relation, Definition \ref{defn:poset}. \textqed

We note that, if $\mu, \nu \in \mathsf{B}(G)$ are distinct, then $(\spt \mu, M_\mu)$ and $(\spt \nu, M_\nu)$ are incomparable with respect to the $\leq$ relation: Lemma \ref{lem:basic-iff} precludes the possibility of strict inequality, and equality is impossible by definition of a basic measure. By this fact and (the contrapositive of) the previous lemma, we conclude as follows.

\begin{cor} \label{cor:not-in-convex-hull}
    Let $S \subseteq \mathsf{B}(G)$ be a compatible set. If $\mu$ is a basic balanced measure not in $S$, then $\mu$ is not in the convex hull of $S$.
\end{cor}

Denote by $\mathcal{P}$ be the power set operator and $C$ the convex hull operator on subsets of $\R^n$. If $\mathcal{S} \subseteq \mathcal{P}(\mathsf{B}(G))$ is the collection of all compatible sets of basic balanced measures on $G$, then Theorem \ref{thm:decomp} tells us that
\begin{equation*}
    \bigcup_{S \in \mathcal{S}} C(S)
\end{equation*}
is precisely the collection $\bbb(G)$ of all balanced measures on $G$. This notation allows us to recast Corollary \ref{cor:not-in-convex-hull} as the aforementioned ``weak minimality" of $\mathsf{B}(G)$.

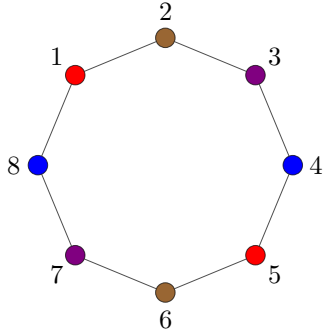
\begin{figure}[h]
    \centering
    \begin{tikzpicture}
        \node at (0,0) {\includegraphics[scale=0.3]{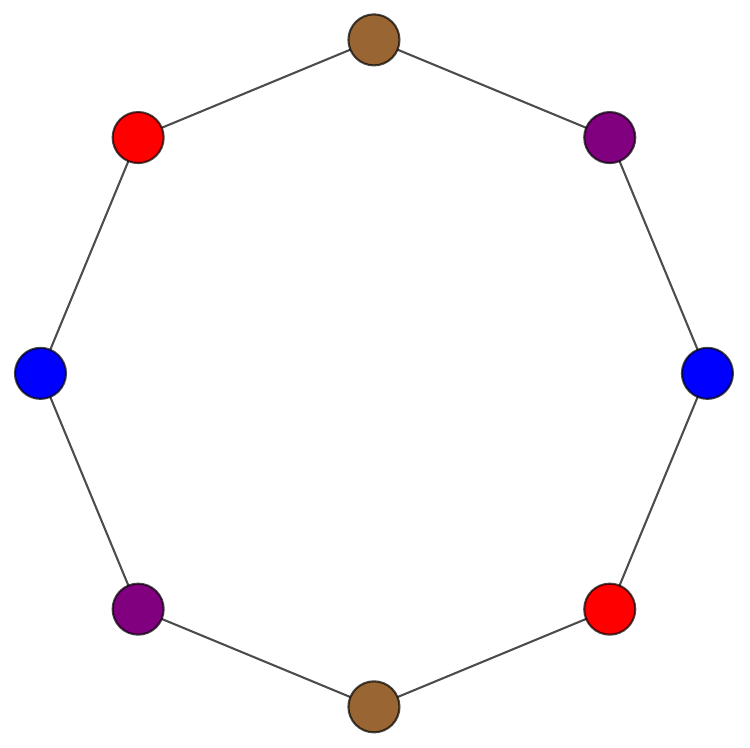}};
        \node at (-1.45,1.45) {$1$};
        \node at (0,2.05) {$2$};
        \node at (1.45,1.45) {$3$};
        \node at (2,0) {$4$};
        \node at (1.45,-1.45) {$5$};
        \node at (0,-2.05) {$6$};
        \node at (-1.45,-1.45) {$7$};
        \node at (-2.02,0) {$8$};
    \end{tikzpicture}
    \caption{The $C_8$ graph equipped with the measures $\mu$, $\nu$, and $\rho$ of Equation \eqref{eq:basic-measures-C4}.}
    \label{fig:C8}
\end{figure}

We illustrate a simple example for Corollary \ref{cor:not-in-convex-hull} below. Let $S = \{\mu, \nu, \rho\}$, where
\begin{equation} \label{eq:basic-measures-C4}
\begin{aligned}
    \mu &= \left( \tfrac{1}{2}, 0, 0, 0, \tfrac{1}{2}, 0, 0, 0 \right) \ \text{(in red)}, \\
    \nu &= \left( 0, 0, \tfrac{1}{2}, 0, 0, 0, \tfrac{1}{2}, 0 \right) \ \text{(in purple)}, \ \text{ and } \\
    \rho &= \left( 0, 0, 0, \tfrac{1}{2}, 0, 0, 0, \tfrac{1}{2} \right) \ \text{(in blue)}
\end{aligned}
\end{equation}
are three compatible basic balanced measures on the $C_8$ graph (Figure \ref{fig:C8}). Since the basic balanced measure 
\begin{equation*}
    \sigma = \left( 0, \tfrac{1}{2}, 0, 0, 0, \tfrac{1}{2}, 0, 0 \right) \ \text{(in brown)}
\end{equation*}
is not in $S$, by Corollary \ref{cor:not-in-convex-hull}, it can not be contained in the convex hull of $S$.

\begin{cor} \label{cor:weak-minimality}
    Let $\mathcal{P}$ and $C$ be as above. Then, for all strict subsets $Q \subsetneq \mathsf{B}(G)$, there is no subcollection $\mathcal{T} \subseteq \mathcal{P}(Q)$ such that $\bigcup_{S \in \mathcal{T}} C(S) = \bbb(G)$.
\end{cor}

\pf Let $Q \subsetneq \mathsf{B}(G)$ and suppose for contradiction that $\bbb(G) = \bigcup_{S \in \mathcal{T}} C(S)$ for some subcollection $\mathcal{T} \subseteq \mathcal{P}(Q)$. We see that $\mathcal{T}$ must only contain compatible sets, as by Corollary \ref{cor:compatible-iff}, the convex hull of any non-compatible set contains measures that are not balanced. Let $\mu \in \mathsf{B}(G) \setminus Q$. Then each $S \in \mathcal{T}$ is a compatible set of basic balanced measures that does not contain $\mu$, so by Corollary \ref{cor:not-in-convex-hull}, $\mu \not\in C(S)$ and, in turn, $\mu \not\in \bigcup_{S \in \mathcal{T}} C(S)$. This is a contradiction. \textqed

Note that we do \textit{not} claim that no element of $\mathsf{B}(G)$ is a convex combination of any others\textemdash hence the ``weakness" of the statement. Counterexamples in which $\mathsf{B}(G) = \bbb(G)$ are fairly abundant.

\section{Cardinality Estimates and Balanced Measures on Graph Joins} \label{s:joins}

\begin{figure}[h!]
\begin{tikzpicture}
        \node at (0,0) {\includegraphics[width=3cm]{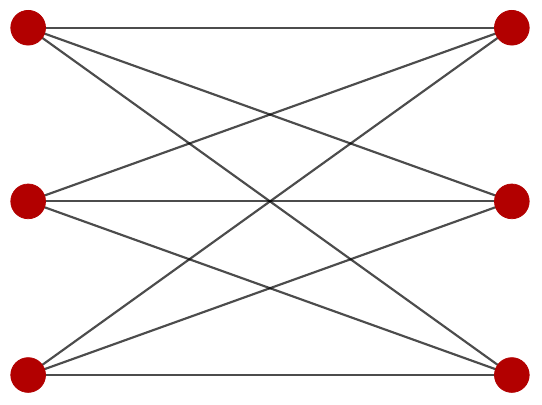}};
        \node at (-3,-3) {\includegraphics[width = 3cm]{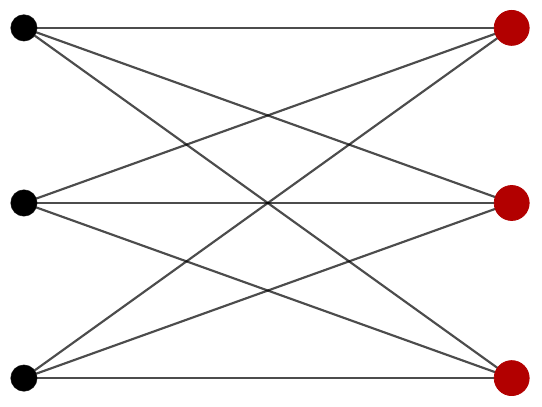}};
        \node at (3,-3) {\includegraphics[width = 3cm]{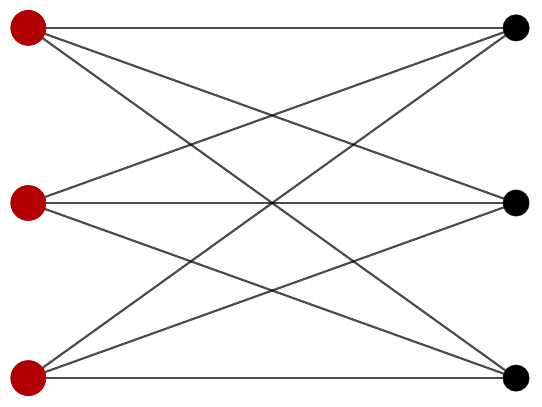}};
    \end{tikzpicture}
    \caption{The bipartite graph $G = \bar{K_3} + \bar{K_3}$ on three vertices. The uniform measure depicted above is a basic balanced measure that can be decomposed into the two basic balanced measures depicted below it. Note that this does not contradict Corollary \ref{cor:weak-minimality}.}
    \label{fig:basic-nonminimal-bipartite}
\end{figure}

In this section we discuss the properties of balanced measures on graph joins (see Figure \ref{fig:basic-nonminimal-bipartite}), which allows us to construct an exponential lower bound on the number of basic balanced measures on a graph. 

\begin{defn}\label{defn:graphjoin}
    Let $G_1$ and $G_2$ be two graphs. The \define{graph join} of $G_1$ and $G_2$ is the graph $G_1 + G_2$ containing all the vertices and edges from $G_1$ and $G_2$, and adding an edge from each vertex of $G_1$ to each vertex of $G_2$.
\end{defn}

Observe that it is natural to consider graph joins when studying distances on graphs due to the following. 

\begin{lem} \label{lem:join-dist}
    Let $G$ be a graph of diameter $2$ and $H$ any other graph. Then we have $d_G(u,v) = d_{G+H}(u,v)$ for all $u,v \in G$, where $d_G$ is the metric on $G$ and $d_{G+H}$ the metric on $G+H$. 
\end{lem}

\pf The graph $G$ has diameter $2$ so $d_G(u,v)=1$ if and only if $u \sim v$ and $d_G(u,v)=2$ otherwise. Likewise, $G+H$ has diameter $2$ so $d_{G+H}$ follows the same rule. Noting that the graph join does not affect the adjacency of elements in $G$, we get the desired result. \textqed

\begin{cor} \label{cor:meas-on-im}
    With $G$ and $H$ as above, there exists a bijection between the balanced measures on $G$ and balanced measures supported in the image of $G$ in $G+H$. 
\end{cor}

\pf Let $\mu_G$ be a balanced measure on $G$. Define a measure $\mu_{G+H}$ on $G+H$ by
\begin{equation*}
    \mu_{G+H}(v) = \left\{ 
    \begin{array}{cl}
        \mu_G(v) & \text{if } v \in G, \\
        0 & \text{otherwise}.
    \end{array} \right.
\end{equation*} 
Then, for every $v \in V$, $T_{\mu_{G+H}}(v) \in \spt \mu_{G+H}$ implies
\begin{equation*}
    T_{\mu_{G+H}}(v) = T_{{\mu}_G}(v) = \max_{w \in V} T_{{\mu}_G}(w) = \max_{w \in V \times V'} T_{\mu_{G+H}}(w)
\end{equation*}
by Lemma \ref{lem:join-dist}. For $v \in V'$, note that $\mu_{G+H}(v,w)=1$ for all $w \in \spt \mu_{G+H} \subseteq G$ so $T_{\mu_{G+H}}(v) = 1 \leq \max_{w \in G+H} T_{\mu_{G+H}}(w)$ hence, $\mu_{G+H}$ is a balanced measure.

Now, given a measure $\mu_{G+H}$ supported on the image of $G$ in $G+H$, define a measure $\mu_G$ on $G$ by $\mu_G(v)=\mu_{G+H}(v)$. By Lemma \ref{lem:join-dist}, we have $T_{\mu_G} = T_{\mu_{G+H}}$ so $\mu_G$ is balanced. Clearly, these constructions give us a bijection so we are done. \textqed

With this in mind, we start with two general lemmas. 

\begin{lem} \label{lem:same-dist-same-mu}
    Let $\mu$ be a balanced measure on a graph $G$. Then, if $v,w$ satisfy the property that for all $u \notin \{v,w\}$, we have $d(u,v)=d(u,w)$, then $\mu(v)=\mu(w)$.
\end{lem}

\pf If $v, w \notin \spt \mu$, then it is trivial so without loss of generality, assume that $v \in \spt \mu$. Then
\begin{equation*}
    0 \leq T_{\mu}(v) - T_{\mu}(w) = d(v,w) \+ \big( \mu(w) - \mu(v) \big),
\end{equation*}
which forces $\mu(w) \geq \mu(v)$. In particular, $\mu(w) > 0$, so $T_{\mu}(v) - T_{\mu}(w) = 0$ and, hence, $\mu(w) = \mu(v)$. \textqed

\begin{lem} \label{lem:conn-to-all}
    Let $\mu$ be a balanced measure on a graph $G$ with $|G| \geq 3$ and let $v$ be a vertex in $G$ such that $v$ is connected to every other vertex. Then, $v \in \spt \mu$ if and only if $G$ is complete. 
\end{lem}

\pf First note that if $G$ is complete, then Lemma \ref{lem:same-dist-same-mu} implies that the only balanced measure on $G$ is the the uniform measure, which is supported everywhere. 

Now let $v \in \spt \mu$ be a vertex connected to every other vertex in $V$, and assume for a contradiction that $G$ is not complete. Then there exist $u,w \in V$ such that $d(u,w) \geq 2$. Observe that
\begin{align*}
  0 &\geq T_{\mu}(u) - T_{\mu}(v) = \sum_{z \in V} \big( d(u,z) - d(v,z) \big) \mu(z) \\
  &= \sum_{z \in V} \big( d(u,z) - 1 \big) \mu(z) + \mu(v) \geq \mu(w) - \mu(u) + \mu(v).  
\end{align*}
So, $\mu(u) \geq \mu(w)+\mu(v)$ and likewise, $\mu(w) \geq \mu(u) + \mu(v)$. However, combining these, we get
\begin{equation*}
    \mu(u) \geq \mu(w) + \mu(v) \geq \mu(u) + 2\mu(v),
\end{equation*}
the desired contradiction. Hence, $G$ must be complete. \textqed

From these lemmas we get the following corollary, which shows that taking the graph join with the singleton graph is a trivial operation for graphs of diameter $2$.

In what follows, let $K_p$ denote the complete graph on $p$ vertices and $\bar{G}$ the complement of a graph $G = (V,E)$, i.e., the graph with vertex set $V$ and edge set complementary to $E$.

\begin{cor} \label{cor:plus-one}
    For any graph $G$ of diameter at most $2$, there is a bijection $\bbb(G) \to \bbb(K_1 + G)$ between the balanced measures on $G$ and those on $K_1 + G$.  
\end{cor}

\pf Once again, by Lemma \ref{lem:same-dist-same-mu}, the only measure on a complete graph is the uniform measure. If $G$ is not complete, then by Lemma \ref{lem:conn-to-all}, the vertex $v$ corresponding to the $K_1$ summand cannot be in the support for a balanced measure. So, any balanced measure on $G+H$ is supported in the image of $G$ and by Corollary \ref{cor:meas-on-im}, we are done. \textqed

These results allow us to construct a graph which will give the aforementioned explicit lower bound on the number of basic measures.

\begin{lem} \label{lem:sum_K_3}
    The graph 
    \[G = \sum_{i=1}^l K_1 + \sum_{i=1}^n \bar{K_3}\] 
    admits exactly $2^n - 1$ balanced measures. 
\end{lem}

\pf By Corollary \ref{cor:plus-one}, it suffices to consider the case $l=0$. Let $B=2J-2I$, where $J$ is the $3 \times 3$ all ones matrix and $I$ is the $3 \times 3$ identity matrix. Then, the distance matrix $D_G$ can be written as the block matrix with $B$ on the diagonals and $J$ off the diagonal. Now let $\mathbbm{1}_n$ be the $1 \times n$ vector of all ones. In addition, let $\{ a_i \}_{i=1}^n \subset (0,1)^n$ such that for some subsequence $\{ a_{i_k} \}_{k=1}^m$, we have $a_{i_k}=1/(3m)$ and for $a_i \notin \{a_{i_k}\}_{k=1}^m$, we have $a_i=0$. Then, we claim that the measure $\mu = (a_1, ..., a_n)^T \otimes \mathbbm{1}_3$ is balanced, where $\otimes$ denotes the Kronecker product. First, observe that
\begin{equation*}
    \|\mu\|_1 = \sum_{i=1}^n 3a_i = \sum_{k=1}^m 3a_{i_k} = \sum_{k=1}^m \frac{1}{m} = 1.
\end{equation*}
So, $\mu$ is a probability measure. Furthermore, note that
\begin{equation*}
    D \mu = \begin{pmatrix}
        4 & 3 & \cdots & 3 \\
        3 & 4 & \cdots & 3 \\
        \vdots & \vdots & \ddots & \vdots \\
        3 & 3 & \cdots & 4
    \end{pmatrix} \begin{pmatrix}
        a_1 \\
        a_2 \\
        \vdots \\
        a_n 
    \end{pmatrix} \otimes \mathbbm{1}_3 = \mu + \mathbbm{1}_n,
\end{equation*}
whence $\mu$ is also balanced. There are $2^n - 1$ distinct choices for $\big\{ a_{i_k} \big\}_{k=1}^m$, so it remains only to show that there are no other balanced measures on $G$. 
    
We start by showing that all balanced measures on $G$ must be of the form $\mu = (a_1, ..., a_n)^T \otimes \mathbbm{1}_3$, meaning in other words, that they have to be constant on each copy of $\bar{K_3}$ under the graph join. Observe that for any $\{ v_1, v_2, v_3 \}$ comprising the image of some $\bar{K_3}$ under the graph join, we have
\begin{equation*}
    d(v_1,v_2) = d(v_2,v_3) = d(v_1,v_3) = 2, \text{whereas}
\end{equation*}
\begin{equation*}
    d(v_1,w) = d(v_2,w) = d(v_3,w) = 1 \qquad \text{whenever } w \not\in \{ v_1, v_2, v_3 \}.
\end{equation*}
Thus, for any balanced measure $\mu$, we have by Lemma \ref{lem:same-dist-same-mu} that $\mu(v_1) = \mu(v_2) = \mu(v_3)$ as desired. 
    
All that is left to show is that $a_i, a_j > 0$ implies that $a_i=a_j$, which would force $\mu$ to be one of the $2^{n}-1$ measures we constructed earlier. The calculation of $D\mu$ performed earlier applies here as well, so we have $D\mu = \mu + \mathbbm{1}_n$. As a result, $\mu$ being balanced immediately entails the desired property. \textqed

\begin{rmk} \label{rmk:not-minimal}
    Because there are only finitely many balanced measures, they must belong to distinct connected components in the space of probability measures, namely singletons. Hence, there are no pairwise compatible measures so they are all basic. This gives us a set of basic measures which is not a minimal generating set for $\bbb(G)$ (via convex combinations) because the measures supported on three vertices generate $\bbb(G)$. 
\end{rmk}

An application of these statements provides the lower bound in Proposition \ref{prop:bounds}.

\noindent \define{Proof of Proposition \ref{prop:bounds}.} The upper bound follows from the fact that every connected component is a convex combination of basic balanced measures, of which there are at most $2^{2n} - 1$. For the lower bound, the case $n = 2$ is trivial and for $n > 2$, by Lemma \ref{lem:sum_K_3}, if $n = 3k + l$ for some $k \in \N$ and $l \in \{0,1,2\}$, then the graph $\sum_{i=1}^l K_1 + \sum_{i=1}^k \bar{K_3}$ has $2^k - 1 = 2^{\lfloor n/3 \rfloor} - 1$ balanced measures, each of which is a basic measure. \textqed

\section{Examples} \label{s:examples}

\subsection{Non-minimal basic balanced measures}

\begin{figure}
    \centering
    \includegraphics[scale=0.33]{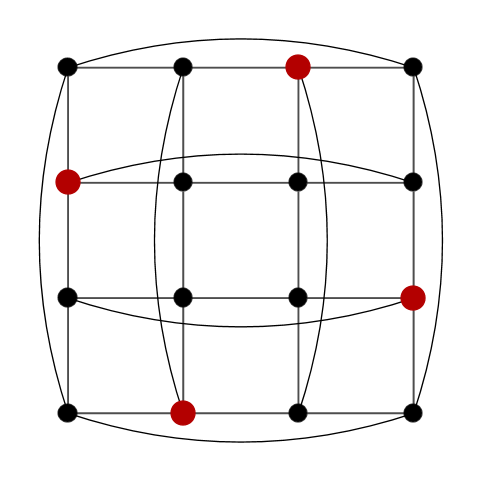}
    \caption{The product graph $C_4 \times C_4$ arranged as a grid and equipped with a balanced measure supported on $4$ vertices\textemdash exactly $1$ from each row and each column.}
    \label{fig:C4xC4}
\end{figure}

Let $C_4 \times C_4$ be the product of the graph $C_4$ with itself, depicted in Figure \ref{fig:C4xC4}. Any measure with weight $1/4$ on each of $4$ distinct vertices\textemdash $1$ from each row and $1$ from each column\textemdash is a basic balanced measure. However, any $2$-vertex subset thereof also supports a basic balanced measure with equal weights. Consequently, a basic balanced measure need not have minimal support among all balanced measures.

\subsection{Limits of Extrapolation}

Lemma \ref{lem:one-sided} shows that given a pair of balanced measures satisfying certain conditions on their support and cost-maximizing set, we may find a bounded line segment through the pair (considered as points in the space of measures on the graph) where every point on the line segment is a balanced measure. In the lemma, we show that all the measures in the interior of the line segment have the same support and cost-maximizing set, and the endpoints both have either smaller support or larger cost-maximizing set. We may conceptualize this as attempting to extrapolate from one point past another, and only stopping when a point is reached where the support decreases in size or the cost-maximizing set increases in size.

With how central Lemma \ref{lem:one-sided} is to the analysis of \S \ref{s:basic-measures}, one might wonder whether both of these ``stopping conditions" are necessary. There are numerous cases where extrapolation must stop because the support decreases in size. For example, extrapolation from $(1/2,0,1/2,0)$ past the uniform measure on $C_4$ will stop at $(0,1/2,0,1/2)$. This exemplifies the support shrinking from $\{1,2,3,4\}$, the uniform measure's support, to $\{2,4\}$.

One might, however, wonder if there are cases where one-sided extrapolation between balanced measures must stop due to an increase in the size of the cost-maximizing set before the support shrinks in size. If this scenario were not possible, we would have a much stronger version of the one-sided extrapolation lemma, with implications on the rest of our results on \S \ref{s:basic-measures}. Unfortunately such scenarios do occur, and we construct such a graph below.

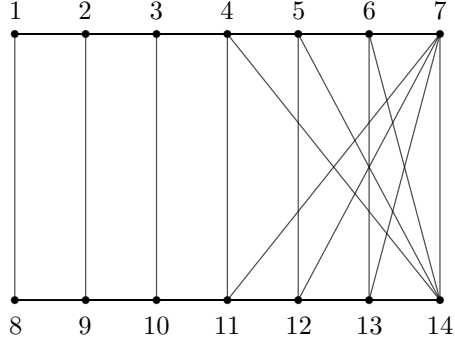
\begin{figure}
    \centering
    \begin{tikzpicture}
        \node at (0,0) {\includegraphics[width=6cm]{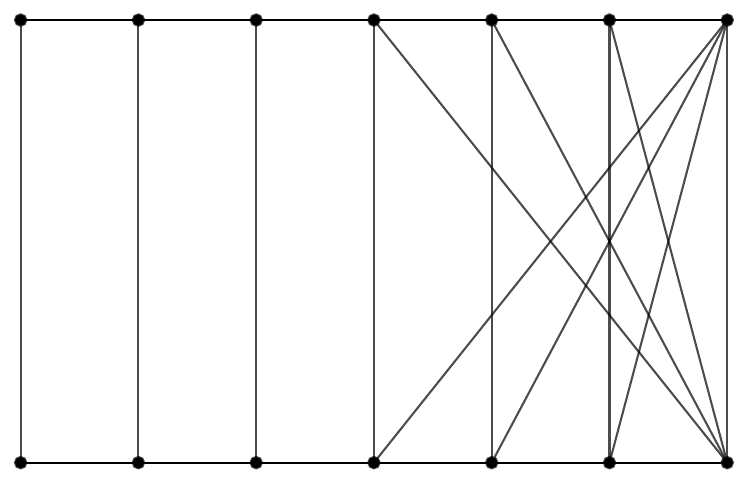}};
        \node at (-2.82,2.1) {$1$};
        \node at (-1.9,2.1) {$2$};
        \node at (-0.95,2.1) {$3$};
        \node at (-0.01,2.1) {$4$};
        \node at (0.94,2.1) {$5$};
        \node at (1.88,2.1) {$6$};
        \node at (2.82,2.1) {$7$};
        \node at (-2.82,-2.1) {$8$};
        \node at (-1.9,-2.1) {$9$};
        \node at (-0.95,-2.1) {$10$};
        \node at (-0.01,-2.1) {$11$};
        \node at (0.94,-2.1) {$12$};
        \node at (1.88,-2.1) {$13$};
        \node at (2.82,-2.1) {$14$};
    \end{tikzpicture}
    \caption{A graph for which one-sided extrapolation stops due to an increase in the size of the cost-maximizing set.}
    \label{fig:balmes}
\end{figure}

Take the graph depicted in Figure \ref{fig:balmes}, where the vertices $\{1, \dots, 7\}$ and $\{8, \dots, 14\}$ form complete graphs (which are not visible due to the choice of embedding due to overlapping edges).

This has the distance matrix
\begin{equation*}
D = \scriptsize
\left(
\begin{array}{cccccccccccccc}
 0 & 1 & 1 & 1 & 1 & 1 & 1 & 2 & 2 & 1 & 2 & 2 & 2 & 2 \\
 1 & 0 & 1 & 1 & 1 & 1 & 1 & 2 & 1 & 2 & 2 & 2 & 2 & 2 \\
 1 & 1 & 0 & 1 & 1 & 1 & 1 & 1 & 2 & 2 & 2 & 2 & 2 & 2 \\
 1 & 1 & 1 & 0 & 1 & 1 & 1 & 2 & 2 & 2 & 2 & 1 & 2 & 1 \\
 1 & 1 & 1 & 1 & 0 & 1 & 1 & 2 & 2 & 2 & 1 & 2 & 2 & 1 \\
 1 & 1 & 1 & 1 & 1 & 0 & 1 & 2 & 2 & 2 & 2 & 2 & 1 & 1 \\
 1 & 1 & 1 & 1 & 1 & 1 & 0 & 2 & 2 & 2 & 1 & 1 & 1 & 1 \\
 2 & 2 & 1 & 2 & 2 & 2 & 2 & 0 & 1 & 1 & 1 & 1 & 1 & 1 \\
 2 & 1 & 2 & 2 & 2 & 2 & 2 & 1 & 0 & 1 & 1 & 1 & 1 & 1 \\
 1 & 2 & 2 & 2 & 2 & 2 & 2 & 1 & 1 & 0 & 1 & 1 & 1 & 1 \\
 2 & 2 & 2 & 2 & 1 & 2 & 1 & 1 & 1 & 1 & 0 & 1 & 1 & 1 \\
 2 & 2 & 2 & 1 & 2 & 2 & 1 & 1 & 1 & 1 & 1 & 0 & 1 & 1 \\
 2 & 2 & 2 & 2 & 2 & 1 & 1 & 1 & 1 & 1 & 1 & 1 & 0 & 1 \\
 2 & 2 & 2 & 1 & 1 & 1 & 1 & 1 & 1 & 1 & 1 & 1 & 1 & 0 \\
\end{array}
\right)\end{equation*}
Then, for $a,b$ such that $6a+6b=1$, consider $\mu_a =(b\mathbbm{1}_3, a \mathbbm{1}_3,0, a \mathbbm{1}_3, b\mathbbm{1}_3,0)^T$
We see
\[D\mu_a=((8a+8b)\mathbbm{1}_6,9a+6b,(8a+8b) \mathbbm{1}_6, 6a+9b)^T.\]
Then $\mu_a$ is balanced if and only if $D\mu_a$ is maximized on the $8a + 8b$ entries. That is, $\mu_a$ is balanced if and only if $8a + 8b \geq 9a + 6b, 6a + 9b$. This occurs when $b/2 < a < 2b$. Since we must have $6a + 6b = 1$, we may conclude through simple computation that $\mu_a$ is balanced if and only if $a \in [1/18,1/9]$.

For $a \in [1/18,1/9]$, we also have $b \in [1/18,1/9]$, so every balanced $\mu_a$ has the same support. However, we may compute explicitly that $D\mu_{1/18}$ is maximized on every vertex except the 8th, $D\mu_{1/9}$ is maximized on every vertex except the 7th, and for $1/18 < a < 1/9$, $D\mu_a$ is maximized on every vertex except both the 7th and 8th.

Thus if we were to extrapolate from $\mu_{1/9}$ past, say, $\mu_{1/12}$ (the measure where $a = b$), we would stop at $\mu_{1/18}$, where the support is the same but the cost-maximizing set has increased in size. This gives an example of the desired type.

\subsection{Basic balanced measures as induced subgraphs}

Consider the following construction.

\begin{defn}
    Given a graph $G$, equip its family $\mathsf{B}(G)$ of basic balanced measures with a graph structure by defining an edge between distinct vertices $\mu,\nu \in \mathsf{B}(G)$ if $\mu$ and $\nu$ are compatible (cf.~Definition \ref{defn:compatible}). We call $\mathsf{B}(G)$ with this graph structure a \define{compatibility graph}.
\end{defn}

Our purpose in this subsection is to show that, for every graph $H$, there exists a graph $G$ such that $\mathsf{B}(G)$ contains an induced subgraph isomorphic to $H$. In this sense, a space of balanced measures can take ``any" geometric configuration, although finding a given configuration may require a large graph $G$.

Recall that a \textit{clique} in a graph is a subset of the vertices that forms a complete induced subgraph, i.e., such that all pairs of vertices in the clique are connected by an edge. We start with a simple lemma that implies a characterization of cliques in $\mathsf{B}(G)$.

\begin{lem} \label{lem:pairs}
    A collection of balanced measures $\mu_1, \dots, \mu_n$ is pairwise compatible if and only if $\{ \mu_1, \dots, \mu_n \}$ is a compatible set of measures. 
\end{lem}

\pf Let $\mu_1, \dots, \mu_n$ be balanced measures that are pairwise compatible. Then, for all $i,j$, we have 
\begin{equation*}
    \spt \mu_i \subseteq \spt \mu_i \cup \spt \mu_j \subseteq M_{\mu_i} \cap M_{\mu_j} \subseteq M_{\mu_j}.
\end{equation*}
It directly follows that
\begin{equation*}
    \bigcup_{i=1}^n \spt \mu_i \subseteq \bigcap_{j=1}^n M_{\mu_j},
\end{equation*}
so $\{ \mu_1, \dots, \mu_n \}$ is a compatible set of measures. 

Conversely, if $\{\mu_1, \dots, \mu_n \}$ is a compatible set of measures, then we have for any $i,j$ that
\begin{equation*}
    \spt \mu_i \cup \spt \mu_j \subseteq \bigcup_{i=1}^n \spt \mu_i \subseteq \bigcap_{i=1}^n M_{\mu_i} \subseteq M_{\mu_i} \cap M_{\mu_j},
\end{equation*}
so $\mu_i$ and $\mu_j$ are compatible. \textqed

\begin{cor}
    A set $\{ \mu_1, \dots, \mu_n \}$ of vertices in $\mathsf{B}(G)$ forms a clique if and only if the corresponding measures are compatible. 
\end{cor}

\pf If the set of vertices forms a clique, then $\mu_1, \dots, \mu_n$ are pairwise compatible, so $\{ \mu_1, \dots, \mu_n \}$ is compatible by Lemma \ref{lem:pairs}. Conversely, if a set $\{ \mu_1, \dots, \mu_n \}$ of basic balanced measures is compatible, then the measures are pairwise compatible. Hence, each pair of measures in $\{ \mu_1, ..., \mu_n \}$ is joined by an edge in $\mathsf{B}(G)$, so they form a clique. \textqed

Given a graph $H$ on $n$ vertices, the graph on $n^2$ vertices defined below will contain a subgraph isomorphic to $H$.

\begin{defn}
    Let $J_n$ be the $n \times n$ matrix of all ones and $I_n$ the $n \times n$ identity matrix. Given a graph $H$ on $n$ vertices and adjacency matrix $A$, construct an $n^2 \times n^2$ block matrix $D$ by replacing each entry $a_{i,j}$ in $A$ with the matrix
    \begin{equation*}
        M_{i,j} \, := \,
        \left\{ \begin{array}{cl}
            2J_3 - 2I_3 & \text{if } i=j \\
            J_3 & \text{if } i \neq j \text{ and } a_{i,j} = 0 \\
            J_3 + I_3 & \text{if } i \neq j \text{ and } a_{i,j} = 1. \\
        \end{array} \right.
    \end{equation*}
    Then define $G_H$ to be the graph with adjacency matrix $2J_{n^2} - 2I_{n^2} - D$.
\end{defn}

\begin{rmk}
    Since $G_H$ has diameter $2$, its distance matrix is precisely $D$.  
\end{rmk}

\begin{lem} \label{lem:vertex-mu}
    For each vertex $v_i \in H$ corresponding to the $i$th row of $A$, let \[\mu_{v_i} = (x_1, \dots, x_n)^T \otimes \mathbbm{1}_3,\] where $x_j = 0$ for all $j \neq i$ and $x_i = 1/3$. Then, $\mu_{v_i}$ is a basic measure on $G_H$. Furthermore, $\mu_{v_i}$ satisfies
    \begin{equation*}
        M_{\mu_{v_i}} = \spt \mu_{v_i} \cup \left( \bigcup_{a_{i,j}=1} \spt \mu_{v_j} \right).
    \end{equation*}
\end{lem}

\pf Define $\mu_{v_i}^j = x_j \mathbbm{1}_3$ for each $i,j$. Observe that if $j \neq i$ then $\mu_{v_i}^j$ is the zero vector. Then, we have
\begin{equation*}
    D\mu = \begin{pmatrix}
    M_{1,1} \+ \mu_{v_i}^1 + \dots + M_{1,n} \+ \mu_{v_i}^n \\
    \vdots \\
    M_{n,1} \mu_{v_i}^1 + \dots + M_{n,n} \+ \mu_{v_i}^n
    \end{pmatrix} = \begin{pmatrix}
    M_{1,i} \+ \mu_{v_i}^1 \\
    \vdots \\
    M_{n,i} \+ \mu_{v_i}^1
    \end{pmatrix}.
\end{equation*} 
Note that
\begin{equation*}
    M_{j,i} \+ \mu_{v_i}^j = \begin{cases}
        4x_i \mathbbm{1}_3 & i=j \\
        3x_i \mathbbm{1}_3 & i \neq j \text{ and } a_{i,j} = 0 \\
        4x_i \mathbbm{1}_3 & i \neq j \text{ and } a_{i,j} = 1, \\
    \end{cases}
\end{equation*}
which gives us that $\mu_{v_i}$ is balanced and $M_{\mu_{v_i}}$ is as desired. Furthermore, it is basic because it has unique support.  \textqed

\begin{prop} \label{prop:subgraph}
    For any graph $H$, there exists a graph $G$ such that $\mathsf{B}(G)$ has an induced subgraph isomorphic to $H$. 
\end{prop}

\pf Let $G = G_H$. Consider the subgraph of $\mathsf{B}(G)$ with vertices $\mu_{v_i}$. By Lemma \ref{lem:vertex-mu}, we know that for $i \neq j$, we have $a_{i,j} = 1$ if and only if $\spt \mu_{v_i} \subseteq M_{\mu_{v_j}}$, which in turn holds if and only if $\{ \mu_{v_i}, \mu_{v_j} \}$ is compatible. Thus, the induced subgraph given by $\{ \mu_{v_1}, \dots, \mu_{v_n} \}$ is clearly isomorphic to $H$ via the map sending $v_i$ to $\mu_{v_i}$. \textqed

\phantomsection
\section*{Acknowledgement}
This paper grew out of a research project at the Washington Experimental Mathematics Laboratory (WXML) at the University of Washington. Our thanks goes to Stefan Steinerberger for his dedicated mentorship throughout the process.

\bibliographystyle{plain}
\bibliography{references}

\end{document}